\documentclass[10pt]{article}
\ifx\pdfoutput\undefined
\usepackage[dvips]{graphicx}
\else
\usepackage[pdftex]{graphicx}
\fi
\catcode`\@=11
\newtheorem{theorem}{Theorem}
\newtheorem{lemma}{Lemma}
\newtheorem{proposition}{Proposition}
\newtheorem{definition}{Definition}
\newtheorem{example}{Example}
\newtheorem{remark}{Remark}

\newtheorem{corollary}{Corollary}
\def\demo{\noindent{\bf Proof .-}}
\def\section{\@startsection {section}{1}{\z@}{-3.5ex plus -1ex
minus-.2ex}{2.3ex plus .2ex}{\normalsize\bf}}

\begin{document}
\ifx\pdfoutput\undefined
\else
\DeclareGraphicsExtensions{.pdf,.gif,.jpg} % the formats we have images in
\fi
\begin{center}
{\Large\bf \textsc{Arithmetical ranks of Stanley-Reisner ideals of simplicial complexes with a cone}}
\end{center}
\vskip.5truecm
\begin{center}
{Margherita Barile\footnote{Partially supported by the Italian Ministry of Education, University and Research.}\\ Dipartimento di Matematica,\\ Universit\`{a} di Bari ``Aldo Moro", Via E. Orabona 4,\\70125 Bari, Italy\footnote{e-mail: barile@dm.uniba.it}}
\end{center}
\begin{center}
{Naoki Terai\\ Department of Mathematics, Faculty of Culture and Education, Saga University,\\ Saga 840-8502, Japan}\footnote {e-mail: terai@cc.saga-u.ac.jp}
\end{center}
\vskip1truecm
\noindent
{\bf Abstract} When a cone is added to a simplicial complex $\Delta$ over one of its faces, we investigate the relation between the arithmetical ranks of the Stanley-Reisner ideals of the original simplicial complex and the new simplicial complex $\Delta'$. In particular, we show that the arithmetical rank of the Stanley-Reisner ideal of $\Delta'$ equals the projective dimension of the Stanley-Reisner ring of $\Delta'$ if the corresponding equality holds for $\Delta$. 

\section*{Introduction and Preliminaries}
The {\it arithmetical rank} (ara) of an ideal $I$ in a commutative Noetherian ring $R$ is the minimal number $s$ of elements $a_1,\dots, a_s$ of $R$ such that $\sqrt I=\sqrt{(a_1,\dots, a_s)}$; one can express this equality by saying that $a_1,\dots, a_s$  {\it generate} $I$ {\it up to radical}. In general height\,$I\leq\,$ara\,$I$; if equality holds, $I$ is called a {\it set-theoretic complete intersection}. 
A better lower bound is provided, in general, by the local cohomological dimension, which, for any squarefree monomial ideal, coincides with the projective dimension (pd), i.e., with the length of any minimal free resolution of the corresponding quotient ring. In this paper we compare this invariant with the arithmetical rank for a certain class of ideals generated by squarefree monomials. 
Let $X$ be a  finite set of indeterminates over a field $K$. A {\it simplicial complex} on $X$ is a set $\Delta$ of subsets of $X$ such that for all $x\in X$, $\{x\}\in\Delta$ and whenever $F\in\Delta$ and $G\subset F$, then $G\in\Delta$. The elements of $\Delta$ are called  {\it faces}, whereas $X$ is called the {\it vertex set} of $\Delta$, and the elements of $X$ are called the {\it vertices} of $\Delta$. If $\Delta$ consists of all subsets of its vertex set, then it is called a {\it simplex}. The simplicial complex $\Delta$ can be associated with an ideal $I_{\Delta}$ of the polynomial ring $R=K[X]$, which is generated by all monomials whose support is not a face of $\Delta$; $I_{\Delta}$ is called the {\it  Stanley-Reisner ideal} of $\Delta$ (over $K$). Its minimal monomial generators are the products of the elements of the minimal non-faces of $\Delta$, and these are squarefree monomials. In fact, this construction provides a one-to-one correspondence between the simplicial complexes on $X$ and the squarefree monomial ideals of $K[X]$ that do not contain  elements of degree one. The quotient ring $K[\Delta]=K[X]/I_{\Delta}$ is called the {\it  Stanley-Reisner ring} of $\Delta$ (over $K$).\newline

In this paper, starting from a simplicial complex $\Delta$, we consider the simplicial complex $\Delta'$ which is the union of $\Delta$ and the simplex ({\it cone}) spanned by a face of $\Delta$ and a new vertex. According to our main theorem,  whenever ara\,$I_{\Delta}=$\,pd\,$K[\Delta]$, then  ara\,$I_{\Delta'}=$\,pd\,$K[\Delta']$. 
Applying this result successively, we inductively deduce that the equality between the arithmetical rank and the projective dimension holds if the Stanley-Reisner ring has a 2-linear resolution, a result which was proven independently, and with completely different arguments, by Morales \cite{Mo}. Our approach is entirely combinatorial and exploits, in the computational part, the linear algebraic techniques developed in \cite{B0}.\newline 
We briefly recall some basic facts about Stanley-Reisner ideals, for which we refer to the extensive treatment given in \cite{BH}, Section 5.  \par\smallskip\noindent 
A maximal face of $\Delta$ is called a {\it facet} of $\Delta$. 
The minimal primes of $I_{\Delta}$ are the ideals of the form
$$P_F=(X\setminus F),\qquad\mbox{ where $F$ is any facet of $\Delta$}.$$
\noindent
It follows that the height of $I_{\Delta}$ (or codimension of $K[\Delta]$) is equal to $\vert X\vert-\max_{F\in\Delta}\vert F\vert$. The number $d=\max_{F\in\Delta}\vert F\vert-1$ is called the {\it dimension} of  $\Delta$; it is evidently equal to $\dim\, K[\Delta]-1$. The ideal $I_{\Delta}$ is unmixed if and only if all facets of $\Delta$ have the same cardinality: we then say that $\Delta$ is {\it pure}. When $\Delta$ is pure, $F\in\Delta$ is called a {\it subfacet} of $\Delta$ if $\vert F \vert = \dim\Delta$.  
If $I_{\Delta}$ is a set-theoretic complete intersection, then one can prove that it is Cohen-Macaulay, from which one concludes that it is unmixed, i.e.,  $\Delta$ is pure. It also follows that the arithmetical rank of $I_{\Delta}$ equals the projective dimension of $K[\Delta]$.
Other classes of squarefree monomial ideals with this property were presented in \cite{B6}, \cite{B7} and \cite{KTY}.\par\smallskip\noindent  
\section{On arithmetical ranks}
We will consider the following two sets of vertices/indeterminates over $K$:  $X=\{x_1,\dots, x_n\}$ and $X'=X\cup\{x_0\}$. We set $R=K[X]$ and $R'=K[X']$. Given the elements $r_1,\dots, r_m\in R$, the ideals they generate in $R$ and $R'$ will be indicated by $(r_1,\dots, r_m)$ and $(r_1,\dots, r_m)R'$ respectively.
\begin{definition}\label{def1}{\rm  Let $F$ be any subset of $X$. We will call {\it cone from $x_0$ over F}, denoted co\,$_{x_0}F$, the simplex on the vertex set $F\cup\{x_0\}$. 
}
\end{definition}
\par\smallskip\noindent
Let $\Delta$ be a simplicial complex on the vertex set $X$, and let $F$ be any face of $\Delta$. Then $\Delta'=\Delta\cup{\rm co}\,_{x_0}F$ is a simplicial complex on the vertex set $X'$. 
\begin{proposition}\label{prop1} If $F=X$, then $\Delta$ and $\Delta'$ are simplices. Otherwise the facets of $\Delta'$ are the facets $G\ne F$ of $\Delta$ and $F\cup\{x_0\}$. Correspondingly, the minimal primes of $I_{\Delta'}$ are the ideals $(P_G+(x_0))R'$ and $P_{F}R'$. 
\end{proposition}
\begin{theorem}\label{theorem1} Let $K$ be algebraically closed. Then, for any face $F$ of $\Delta$,
$${\rm ara}\,I_{\Delta'}\leq \max({\rm ara}\,I_{\Delta}+1, n-\vert F\vert).$$
If $F\ne X$ and $I_{\Delta}$ is a set-theoretic complete intersection, then equality holds.
\end{theorem} 
\demo 
If $F=X$, then the ideals $I_{\Delta}$ and $I_{\Delta'}$ are both zero, so that the claim is trivially true. So assume that $F\ne X$.
Let $G$ be a facet of $\Delta$ containing $F$. Up to renaming the indeterminates, we may assume that $G=\{x_{s+1}, \dots, x_n\}$ and $F=\{x_{t+1}, \dots, x_n\}$ for some integers $s$ and $t$ such that $1\leq s\leq t< n$. Then $P_G=(x_1,\dots, x_s)$ is a minimal prime of $I_{\Delta}$, and $t=n-\vert F\vert$. Moreover,
\begin{equation}\label{R'} I_{\Delta'}=(I_{\Delta}+(x_0x_1,\dots, x_0x_t))R',
\end{equation}
and, by Proposition \ref{prop1},
\begin{equation}\label{PF} I_{\Delta'}\subset(x_1,\dots, x_t)R'.\end{equation}
\noindent
 Set $h=\,{\rm ara}\,I_{\Delta}$, and let $q_1,\dots, q_h\in R$ be polynomials generating $I_{\Delta}$ up to radical. 
For all indices $i\in\{1,\dots, h\}$, we have that $q_i\in I_{\Delta}$, whence $q_i\in (x_1,\dots, x_s)$, i.e., $q_i=\sum_{j=1}^s a_{ij}x_j$ for some $a_{ij}\in R$. Recall that a polynomial belongs to a given monomial ideal   if and only if each of its monomial terms is  divisible by some monomial generator of this ideal. Therefore, up to eliminating redundant terms, we may assume that 
\begin{equation}\label{0}\qquad\qquad a_{ij}x_j\in I_{\Delta}\qquad\qquad\mbox{ for all indices } i\in\{1,\dots, h\},\ j\in\{1,\dots, s\}.\end{equation}
\noindent
Consider the following ring homomorphism:
$$\phi:R\longrightarrow R$$
$$f(x_1,\dots, x_n)\mapsto f(x_1^2,\dots, x_n^2),$$
\noindent
and, for all $i=1,\dots, h$, set $\bar q_i=\phi(q_i)$.
Moreover, for all $i\in\{1,\dots, h\}$ and all $j\in\{1,\dots, s\}$, set $\bar a_{ij}=\phi(a_{ij})x_j$, so that, for all $i\in\{1,\dots, h\}$, 
\begin{equation}\label{barq}\bar q_i=\sum_{j=1}^s\bar a_{ij}x_j.\end{equation}
\noindent
The monomial terms of every $\bar a_{ij}$ are all of the form $b^2x_j$, where $b$ is a monomial term of $a_{ij}$; hence, by (\ref{0}), $bx_j\in I_{\Delta}$, so that $b^2x_j\in I_{\Delta}$, and, for all $i\in\{1,\dots, h\}$, $j\in\{1,\dots, s\}$, 
\begin{equation}\label{2}\bar a_{ij}\in I_{\Delta},\qquad\mbox{whence}\qquad \bar q_i\in I_{\Delta}.\end{equation}
\noindent
Let $\bar J=(\bar q_1,\dots, \bar q_h).$ We show that 
\begin{equation}\label{1} \sqrt{\bar J}=I_{\Delta}.\end{equation}
\noindent
In view of (\ref{2}) we only need to show inclusion $\supset$. Let $g$ be any monomial generator of $I_{\Delta}$. Since, by assumption, $\sqrt{(q_1,\dots, q_h)}=I_{\Delta}$, for some positive integer $e$ we have that $g^e\in (q_1,\dots, q_h)$, i.e., there are $f_1,\dots, f_h\in R$ such that $g^e=\sum_{i=1}^hf_iq_i$. But then 
$$g^{2e}=\phi(g^e)=\sum_{i=1}^h\phi(f_i)\bar q_i,$$
\noindent
which shows that $g\in\sqrt{\bar J}$. Thus $I_{\Delta}\subset\sqrt{\bar J}$, which completes the proof of (\ref{1}).\newline
For the remaining part of the proof we need to distinguish between the following cases.
\par\smallskip\noindent
\underline{Case 1}. Suppose that $h+1>t.$ We show that ara\,$I_{\Delta'}\leq h+1$.  Consider the $t\times t$ matrix $\bar A=(\bar a_{ij})_{i,j=1,\dots, t}$, where $\bar a_{ij}=0$ whenever $j>s$.
 Let $A_1=\bar A + x_0Id_t$, where $Id_t$ denotes the $t\times t$ identity matrix. 
Set
$$J_1=( \det A_1-x_0^t,\bar q_1+x_0x_1, \dots, \bar q_t+x_0x_t, \bar q_{t+1},\dots, \bar q_h)R'.$$
\noindent
We prove that
\begin{equation}\label{claim11} I_{\Delta'}=\sqrt{J_1}.
\end{equation}
\noindent
Note that  $\det A_1-x_0^t$ is the sum of products each of which involves at least one entry of $\bar A$ as a factor; in view of (\ref{2}), it follows that  $\det A_1-x_0^t\in I_{\Delta}$. This, together with (\ref{R'}) and (\ref{2}), allows us to conclude that $I_{\Delta'}\supset\sqrt{J_1}$. We show the opposite inclusion by means of Hilbert's Nullstellensatz. Let ${\bf x}=(x_0, x_1,\dots, x_n)\in K^{n+1}$ be such that all elements of $J_1$ vanish at ${\bf x}$. We show that ${\bf x}$ annihilates all elements of $I_{\Delta'}$. In the rest of the proof, we shall identify each polynomial with its value at ${\bf x}$. Our assumption can be formulated in the form:
\begin{eqnarray}\label{equationa} \det A_1-x_0^t&=&0,\\
\bar q_1+x_0x_1=\cdots=\bar q_t+x_0x_t&=&0,\label{equationb}\\
\bar q_{t+1}=\cdots =\bar q_h&=&0.\label{equationc}
\end{eqnarray}
\noindent
 We distinguish between two cases. First suppose that $\det A_1\ne 0$. Note that $A_1$ is the matrix of coefficients of the square system of homogeneous linear equations   
$$\sum_{j=1}^t(\bar a_{ij}+\delta_{ij}x_0)y_j=0\qquad(i=1,\dots, t)$$
in the unknowns $y_1,\dots, y_t$. By Cramer's Rule it only has  the trivial solution. Therefore, in view of (\ref{barq}), (\ref{equationb}) implies that $x_1=\cdots=x_t=0$. In view of (\ref{PF}), it then follows that ${\bf x}$ annihilates all elements of $I_{\Delta'}$. Now suppose that $\det A_1=0$. Then  from (\ref{equationa}) we deduce that $x_0=0$, so that from (\ref{equationb}) and (\ref{equationc}) we further have that $\bar q_1=\cdots=\bar q_h=0$. These equalities, together with (\ref{1}), imply that all elements of $I_{\Delta}$ vanish at ${\bf x}$.  In view of (\ref{R'}), we again conclude that ${\bf x}$ annihilates all elements of $I_{\Delta'}$. This completes the proof of (\ref{claim11}), which shows that in Case 1 ara\,$I_{\Delta'}\leq h+1$, as required. \newline
\underline{Case 2}. Suppose that $h+1\leq t.$ We show that ara\,$I_{\Delta'}\leq t$.  Consider the $h\times h$ matrix $\bar A'=(\bar a_{ij})_{i,j=1,\dots, h}$, where $\bar a_{ij}=0$ whenever $j>s$.\newline
\underline{Case 2.1}. First assume that char\,$K=p>0$. Let $A_{21}=\bar A' + x_0(x_0-x_{h+1})Id_h$. Choose a positive integer $\ell$ such that $p^{\ell}>h$, and set
\begin{eqnarray*}J_{21}=&&(x_0^{p^{\ell}-h}(x_0-x_{h+1})^{p^{\ell}-h}\det A_{21}-x_0^{2p^{\ell}},\\
&&\bar q_1+x_0(x_0-x_{h+1})x_1,\dots,  \bar q_h+x_0(x_0-x_{h+1})x_h,\\
&& x_0x_{h+2},\dots, x_0x_t)R'.
\end{eqnarray*}
\noindent
We prove that
\begin{equation}\label{claim21} I_{\Delta'}=\sqrt{J_{21}}.
\end{equation}
\noindent
Note that $x_0^{p^{\ell}-h}(x_0-x_{h+1})^{p^{\ell}-h}\det A_{21}-x_0^{2p^{\ell}}$ is the sum of products each of which involves at least one entry of $\bar A'$ or the monomial $x_0x_{h+1}$ as a factor. This, together with (\ref{R'}) and (\ref{2}), allows us to conclude that $I_{\Delta'}\supset\sqrt{J_{21}}$. We show the opposite inclusion by means of Hilbert's Nullstellensatz. Let ${\bf x}\in K^{n+1}$ be such that all elements of $J_{21}$ vanish at ${\bf x}$. We show that ${\bf x}$ annihilates all elements of $I_{\Delta'}$.
Our assumption can be formulated in the form:
\begin{eqnarray}\label{equationa''} x_0^{p^{\ell}-h}(x_0-x_{h+1})^{p^{\ell}-h}\det A_{21}-x_0^{2p^{\ell}}&=&0,\\
\bar q_1+x_0(x_0-x_{h+1})x_1=\cdots=\bar q_h+x_0(x_0-x_{h+1})x_h&=&0,\label{equationb''}\\
x_0x_{h+2}=\cdots = x_0x_t&=&0.\label{equationc''}
\end{eqnarray}
\noindent
 We distinguish between two cases. First suppose that $\det A_{21}\ne 0$. Note that $A_{21}$ is the matrix of coefficients of the square system of homogeneous linear equations   
$$\sum_{j=1}^h(\bar a_{ij}+\delta_{ij}x_0(x_0-x_{h+1}))y_j=0\qquad(i=1,\dots, h)$$
in the unknowns $y_1,\dots, y_h$. By Cramer's Rule it only has  the trivial solution. Therefore, in view of (\ref{barq}), (\ref{equationb''}) implies that $x_1=\cdots=x_h=0$, whence $\bar q_1=\cdots=\bar q_h=0$. By (\ref{1}) it follows that ${\bf x}$ annihilates all elements of $I_{\Delta}$. In view of (\ref{R'}) and (\ref{equationc''}), there remains to show that $x_0x_{h+1}=0$. Since $\bar A '=0$, we have that $A_{21}= x_0(x_0-x_{h+1})Id_h$, whence $\det A_{21}=x_0^h(x_0-x_{h+1})^h$, and $x_0^{p^{\ell}-h}(x_0-x_{h+1})^{p^{\ell}-h}\det A_{21}-x_0^{p^{\ell}}(x_0-x_{h+1})^{p^{\ell}}=0$. In view of  (\ref{equationa''}), this implies that $x_0x_{h+1}=0$, as claimed. Now suppose that $\det A_{21}=0$. Then from (\ref{equationa''}) we deduce that $x_0=0$, which, in turn, according to (\ref{equationb''}), gives us $\bar q_1=\cdots=\bar q_h=0$. In view of (\ref{R'}) and (\ref{1}) we have that ${\bf x}$ annihilates all elements of $I_{\Delta'}$. This proves (\ref{claim21}).\newline 
\underline{Case 2.2}. Now assume that char\,$K=0$. 
Let $\omega_1,\dots, \omega_h$ be the $h$-th roots of unity in $K$. Set 
$$A_{22}=\bar A' + x_0(x_0-\omega_1x_{h+1}, \dots, x_0-\omega_hx_{h+1})Id_h,$$
\noindent and
\begin{eqnarray*}
J_{22}=&&(\det A_{22}-x_0^{2h}, \\
&&\bar q_1+x_0(x_0-\omega_1x_{h+1})x_1,\dots, \bar q_h+x_0(x_0-\omega_hx_{h+1})x_h,\\
&&x_0x_{h+2},  \dots, x_0x_t).
\end{eqnarray*}
\noindent
We prove that
\begin{equation}\label{claim22} I_{\Delta'}=\sqrt{J_{22}}.
\end{equation}
\noindent
It suffices to show that $I_{\Delta'}\subset \sqrt{J_{22}}$.   We use Hilbert's Nullstellensatz. Let ${\bf x}\in K^{n+1}$ be such that all elements of $J_{22}$ vanish at ${\bf x}$. We show that ${\bf x}$ annihilates all elements of $I_{\Delta'}$.  Our assumption can be formulated in the form:
\begin{eqnarray}\label{equationa'} \det A_{22}-x_0^{2h}&=&0,\\
\bar q_1+x_0(x_0-\omega_1x_{h+1})x_1=\cdots=\bar q_h+x_0(x_0-\omega_hx_{h+1})x_h&=&0,\label{equationb'}\\
x_0x_{h+2}=\cdots=x_0x_t&=&0.\label{equationc'}
\end{eqnarray}
\noindent
 We distinguish between two cases. First suppose that $\det A_{22}\ne 0$. Note that $A_{22}$ is the matrix of coefficients of the square system of homogeneous linear equations   
$$\sum_{j=1}^h(\bar a_{ij}+\delta_{ij}x_0(x_0-\omega_jx_{h+1}))y_j=0\qquad(i=1,\dots, h)$$
in the unknowns $y_1,\dots, y_h$. By Cramer's Rule it only has  the trivial solution. Therefore (\ref{equationb'}) implies that $x_1=\cdots=x_h=0$, so that $\bar q_1=\cdots=\bar q_h=0$. In view of (\ref{R'}), (\ref{1}) and (\ref{equationc'}),  there remains to show that $x_0x_{h+1}=0$. Again we have that $\bar A'=0$, so that 
$$\det A_{22}= \prod_{i=1}^hx_0(x_0-\omega_ix_{h+1})=x_0^h(x_0^h-x_{h+1}^h).$$
Thus (\ref{equationa'}) implies that $x_0x_{h+1}=0$, as claimed. 
Now suppose that $\det A_{22}=0$. Then from (\ref{equationa'}) we deduce that $x_0=0$, so that from (\ref{equationb'}) we further have that $\bar q_1=\cdots=\bar q_h=0$. These equalities, together with (\ref{R'}) and (\ref{1}), imply that all elements of $I_{\Delta'}$ vanish at ${\bf x}$, which shows our claim (\ref{claim22}).\newline
This proves that in Case 2 ara\,$I_{\Delta'}\leq t$, as required, 
and completes the proof of the first part of the claim. 
Now suppose that $F\ne X$ and that $I_{\Delta}$ is a set-theoretic complete intersection. In this case $h$ is the  height of all minimal primes of $I_{\Delta}$. On the other hand, there is a facet $G'$ of $\Delta$ such that $G'\neq F$. Then, by Proposition \ref{prop1}, $(P_{G'}+(x_0))R'$ is a minimal prime of $I_{\Delta'}$, of height $h+1$. By \cite{M}, Theorem 13.5, we have that ara\,$I_{\Delta'}\geq h+1$. Moreover, by Proposition \ref{prop1}, $P_FR'$ is a minimal prime of $I_{\Delta'}$ of height $n-\vert F\vert$, so that, by  \cite{M}, Theorem 13.5, we also have that ara\,$I_{\Delta'}\geq n-\vert F\vert.$  Thus ara\,$I_{\Delta'}\geq\max({\rm ara}\,I_{\Delta}+1, n-\vert F\vert).$ By the first part of the claim, it follows that equality holds. This completes the proof.
\par\medskip\noindent
\begin{remark}{\rm Case 2.2 of the proof of Theorem \ref{theorem1} applies, more in general, whenever ara\,$I_{\Delta}$ is coprime with respect to char\,$K$.
}
\end{remark}
As a special case of Theorem \ref{theorem1} we get the following result, which was already proven in \cite{B0}, Theorem 1. 
\begin{corollary}\label{coro1}
Suppose that $\Delta$ is not a simplex. If $F$ is a facet of $\Delta$ and $I_{\Delta}$ is a set-theoretic complete intersection, then 
$${\rm ara }\,I_{\Delta'}=\,{\rm ara}\, I_{\Delta}+1.$$
\end{corollary}
\demo Since $\Delta$ is pure, we have that ${\rm ara}\, I_{\Delta}=\,{\rm height}\,I_{\Delta}=n-\vert F\vert$. It suffices to apply the second part of Theorem \ref{theorem1} to conclude.
\par\medskip\noindent
From Theorem 1 we also deduce the following result. 
\begin{corollary}\label{coro2}
Suppose $I_{\Delta}$ is a set-theoretic complete intersection. If $F$ is a subfacet of $\Delta$, then $I_{\Delta'}$ is also a set-theoretic complete intersection. 
\end{corollary}
\demo The ideal $I_{\Delta}$ has pure height equal to ara\,$I_{\Delta}$. Moreover, by Proposition \ref{prop1},  height\,$I_{\Delta'}=$\,height\,$I_{\Delta}+1$. On the other hand, $n-\vert F\vert=$\,height\,$I_{\Delta}+1$, so that, by Theorem \ref{theorem1}, ara\,$I_{\Delta'}\leq$\,ara\,$I_{\Delta}+1$, whence ara\,$I_{\Delta'}\leq$\,height\,$I_{\Delta'}$. But the opposite inequality is always true. This completes the proof. 
\begin{example}\label{example1}{\rm
Consider the simplicial complex $\Delta$ on the vertex set $\{x_1,\dots, x_4\}$ whose facets are $\{x_1, x_2\}$, $\{x_2, x_3\}$, $\{x_3, x_4\}$, $\{x_4, x_1\}$. 
Its Stanley-Reisner ideal in the polynomial ring $R=K[x_1,\dots, x_4]$ is 
$$I_{\Delta}=(x_1x_3,\ x_2x_4), $$
\noindent
whose minimal prime decomposition  is
$$I=(x_3, x_4)\cap(x_1, x_4)\cap(x_1, x_2)\cap(x_2, x_3).$$
\noindent It is a complete intersection.
Let $F=\{x_4\}$. The corresponding simplicial complex $\Delta'=\Delta\cup\,{\rm co}\,_{x_0}F$ on the vertex set $\{x_0, x_1,\dots, x_4\}$ has the following facets: $\{x_0, x_4\}$, $\{x_1, x_2\}$, $\{x_2, x_3\}$, $\{x_3, x_4\}$, $\{x_4, x_1\}$.
The monomial generators of $I_{\Delta'}$ are:
$$x_1x_3, x_2x_4,  x_0x_1, x_0x_2, x_0x_3.$$
\noindent
According to Corollary \ref{coro2}, ara\,$I_{\Delta'}=\,$height\,$I_{\Delta'}=3$; we construct three elements generating $I_{\Delta'}$ up to radical applying the procedure described in the proof of Theorem \ref{theorem1}.
With respect to the notation introduced there, we have that $q_1=x_1x_3$, $q_2=x_2x_4$, $h=2$, $P_{F}=(x_1, x_2, x_3)$,  $t=3$. Hence Case 2 applies. We can take $G=\{x_3, x_4\}$ as a facet of $\Delta$ containing $F$, so that $s=2$. We thus have 
 $$\bar A'=\left(\begin{array}{cc}
x_1^2x_3&0\\
0&x_2^2x_4
\end{array}\right)
$$
\noindent
First assume that char\,$K= 2$. Then we can take $\ell=2$, so that
\begin{eqnarray*}
A_{21}&=&\bar A' +x_0(x_0+x_3)Id_2\\
&=&\left(\begin{array}{cc}
x_1^2x_3+x_0(x_0+x_3)&0\\
0&x_2^2x_4+x_0(x_0+x_3)
\end{array}\right).
\end{eqnarray*}
\noindent
It follows that $I_{\Delta'}$ is generated up to radical by the following three elements:
\begin{eqnarray*}&&x_0^2(x_0+x_3)^2\det A_{21}+x_0^8=
x_0^4x_1^2x_2^2x_3x_4 + x_0^2x_1^2x_2^2x_3^3x_4 + x_0^5x_1^2x_3^2 +\\
&&\qquad\qquad\qquad x_0^4x_1^2x_3^3 + x_0^3x_1^2x_3^4 +x_0^6x_1^2x_3 + x_0^6x_2^2x_4 +x_0^5x_2^2x_3x_4 +\\ 
&&\qquad\qquad\qquad x_0^4x_2^2x_3^2x_4 +x_0^3x_2^2x_3^3x_4+x_0^4x_3^4,\\
&&x_1^2x_3^2+x_0^2x_1+x_0x_1x_3,\\
&&x_2^2x_4^2+x_0^2x_2+x_0x_2x_3.
\end{eqnarray*}
 Now assume that char\,$K\neq 2$. Then
\begin{eqnarray*}
A_{22}&=&\bar A' +(x_0(x_0-x_3), x_0(x_0+x_3))Id_2\\
&=&\left(\begin{array}{cc}
x_1^2x_3+x_0(x_0-x_3)&0\\
0&x_2^2x_4+x_0(x_0+x_3)
\end{array}\right).
\end{eqnarray*}
\noindent
It follows that $I_{\Delta'}$ is generated up to radical by the following three elements:
\begin{eqnarray*}&&\!\!\!\!\!\!\!\det A_{22}-x_0^4=
x_1^2x_2^2x_3x_4+x_0^2x_2^2x_4-x_0x_2^2x_3x_4+x_0^2x_1^2x_3+x_0x_1^2x_3^2-x_0^2x_3^2,\\
&&x_1^2x_3^2+x_0^2x_1-x_0x_1x_3,\\
&&x_2^2x_4^2+x_0^2x_2+x_0x_2x_3.
\end{eqnarray*}
\noindent}
\end{example}
\section{On projective dimensions}
For all $Y\subset X$ (and $Y'\subset X'$) let $\Delta_Y$ ($\Delta'_{Y'}$) denote the subcomplex induced by $\Delta$ on $Y$ (by $\Delta'$ on $Y'$). According to Hochster's formula (see \cite{H}, Theorem 5.1, p. 194), the Betti numbers of the Stanley-Reisner ring of a simplicial complex can be completely characterized in terms of the reduced simplicial homology of the induced subcomplexes of $\Delta$: for all nonnegative indices $i$
\begin{equation}\label{Hochster}
\beta_i(K[\Delta'])=\dim_K\,{\rm Tor}_{i}^{R'}(K[\Delta'],K)=
\sum_{Y'\subset X'}\dim_K\tilde H_{\vert Y'\vert -i-1}(\Delta'_{Y'};K).
\end{equation} 
We will use this formula in the proof of the next result.
\begin{lemma}\label{theorem2}
Suppose that $F\neq X$. Then
$${\rm pd}\,K[\Delta']=\max({\rm pd}\,K[\Delta]+1, n-\vert F\vert).$$
\end{lemma}
\demo According to (\ref{Hochster}), 
\begin{equation}\label{pd}{\rm pd}\,K[\Delta']=\max\{i\mid \tilde H_{\vert Y'\vert -i-1}(\Delta'_{Y'};K)\ne 0\mbox{ for some }Y'\subset X'\}.\end{equation}
\noindent
For the sake of simplicity, in the sequel we will omit the coefficient field $K$ in the homology groups. 
If $Y'\subset X$, then $\Delta'_{Y'}=\Delta_{Y'}$. Otherwise, let $Y=Y'\setminus\{x_0\}$. Then, if $F\cap Y\neq\emptyset$, we have that $\Delta'_{Y'}=\Delta_Y\cup\,{\rm co}_{x_0}(F\cap Y)$ is homotopically equivalent to $\Delta_Y$; if $F\cap Y=\emptyset$, we have that the geometric realization of $\Delta'_{Y'}$ is the disjoint union of that of $\Delta_Y$ and $\{x_0\}$. It follows that, for all indices $j$, 
$$\tilde H_j(\Delta'_{Y'})=\left\{
\begin{array}{ll}
\tilde H_j(\Delta_{Y'}) & \mbox{if } Y'\subset X,\\
\tilde H_j(\Delta_{Y}) & \mbox{if } x_0\in Y', \\
&\mbox{ and }F\cap Y\neq\emptyset,\mbox{ or }F\cap Y=\emptyset\mbox{ and }j>0,\\
\tilde H_0(\Delta_{Y})\oplus K &\mbox{otherwise}. 
\end{array}
\right.
$$
This implies that
$$\tilde H_{\vert Y'\vert-i-1}(\Delta'_{Y'})=\left\{
\begin{array}{ll}
\tilde H_{\vert Y\vert-i-1}(\Delta_{Y}) & \mbox{if } Y'=Y\subset X,\\
\tilde H_0(\Delta_{Y})\oplus K& \mbox{if } x_0\in Y',\\
& \mbox{ and }F\cap Y=\emptyset, i=\vert Y\vert=\vert Y'\vert -1,\\
\tilde H_{\vert Y\vert-(i-1)-1}(\Delta_{Y})  &\mbox{otherwise}. 
\end{array}
\right.
$$
\noindent
Hence (\ref{pd}) can be rewritten as
\begin{eqnarray*}
{\rm pd}\,K[\Delta']&=&\max(\max\{i\mid \tilde H_{\vert Y\vert -i-1}(\Delta_{Y})\ne 0\mbox{ for some }Y\subset X\}+1,\\
&&\qquad\max\{ \vert Y'\vert-1\mid x_0\in Y'\subset X', F\cap Y=\emptyset\})\\
&=&\max\{{\rm pd}\,K[\Delta]+1, \vert X'\setminus F\vert-1\}\\
&=&\max\{{\rm pd}\,K[\Delta]+1, n-\vert F\vert\},
\end{eqnarray*}
as claimed.
\par\medskip\noindent
\begin{theorem}\label{theorem3} Let $K$ be algebraically closed.  If {\rm ara}\,$I_{\Delta}=$\,{\rm pd}\,$K[\Delta]$, then \linebreak {\rm ara}\,$I_{\Delta'}=$\,{\rm pd}\,$K[\Delta']$. 
\end{theorem}
\demo If $F=X$, then both $\Delta$ and  $\Delta'$ are simplices, and 
${\rm ara}\ I_{\Delta'} =0=\,{\rm pd}\ K[\Delta']$. Otherwise
we have that 
\begin{eqnarray*}
{\rm ara}\,I_{\Delta'}&\leq&\max({\rm ara}\,I_{\Delta}+1, n-\vert F\vert)\\
&=&\max({\rm pd}\,K[\Delta]+1, n-\vert F\vert)\\
&=&\,{\rm pd}\,K[\Delta']\leq\,{\rm ara}\,I_{\Delta'},
\end{eqnarray*}
where the first inequality and the second equality follow from Theorem \ref{theorem1} and Lemma \ref{theorem2} respectively.
\section{Applications} 
In this section we present applications of Theorem \ref{theorem1} and Theorem \ref{theorem3}. 
Before stating the first result, we need to recall a definition and a lemma. 
\begin{definition} {\rm A  {\it generalized tree} is defined inductively as follows.\newline
(i) A simplex is a generalized tree. \newline
(ii) If $\Delta$ is a generalized tree, then so is 
$\Delta \cup \mbox{\rm co } _{x_0} F$  
for any $F \in \Delta$ and for any new vertex $x_0$. }
\end{definition}
The next claim is a consequence of \cite{F}, Theorem 1, and \cite{D}, Proposition 5.5.1. 
\begin{lemma}\label{DiestelFroberg} 
For any simplicial complex $\Delta$ which is not a simplex,
the Stanley-Reisner ring $K[\Delta]$ of $\Delta$ 
has a $2$-linear resolution if and only if $\Delta$ is a generalized tree.
\end{lemma}
\par\medskip\noindent
From Theorem \ref {theorem3} we can now deduce the following result by Morales (\cite{Mo}, Theorem 8).
\begin{corollary}\label{coro3} Let $K$ be algebraically closed. If $K[\Delta]$ has a $2$-linear resolution, then 
$${\rm ara}\,I_{\Delta}=\,{\rm pd}\,K[\Delta].$$
\end{corollary}
\demo   According to Lemma 2, we have to prove ${\rm ara}\ I_\Delta =\,{\rm pd}\ K[\Delta]$ for a generalized tree $\Delta$ by
induction on the number $|X|$ of vertices.
If $|X|=1$, then $\Delta$ is a simplex, and ara $I_\Delta =0=\,{\rm pd}\ K[\Delta]$.
Assume $|X| >1$. Then there exist a vertex $x \in X$, a generalized tree $\bar{\Delta}$
on the vertex set $X\setminus \{x\}$, and a face $F \in \bar\Delta$
such that $\Delta =\bar{\Delta} \cup \mbox{\rm co } _{x} F$.
By the induction hypothesis, we have ${\rm ara}\ I_{\bar{\Delta}} =\,{\rm pd}\ K[\bar{\Delta}]$.
Hence by Theorem 2, we have  ${\rm ara}\ I_\Delta = {\rm pd}\ K[\Delta]$. 
Thus the claim is proven by induction on the number of vertices.
\par\medskip\noindent
 Next we apply Theorem \ref{theorem3} to compute the arithmetical rank of the Stanley-Reisner ideal 
$I_{\Delta}$ 
of a pure strongly connected simplicial 
complex $\Delta$   which  satisfies the inequality
${\rm reg } I_{\Delta} = \deg K[\Delta]-{\rm  codim }K[\Delta] +1 \ge 4$. Here deg denotes the multiplicity, and reg denotes the Castelnuovo-Mumford regularity. We recall that it is equal to 2 when the ring has a 2-linear resolution. 
\newline
First we introduce some notation.
We say that a pure complex $\Delta$ is {\it  strongly connected}
if 
for any two facets  $F$ and $G$, there exists a sequence of facets
$$
F=F_{0},F_{1}, \dots , F_{m}=G
$$
such that $ F_{i-1}\cap F_{i}$ is a  subfacet for all $i=1,2, \dots , m.$
\newline
Moreover, we denote by $\Delta (n)$ the {\it elementary $(n-1)$-simplex} $\mbox{\bf  2}^{X}$,
which is the power set of $X$,
and put $\Delta (0)= \{ \emptyset \}$. 
We also put 
$\partial {\Delta (n)}=\Delta (n) \setminus \{ X\}$,
and call this the {\it boundary complex} of $\Delta (n)$.
\newline
Let $\Delta _{0}$  be a pure $d$-dimensional simplicial 
complex.
Take a subfacet $F_0$ in $\Delta _{0}$  and a new vertex $y_1$.
Let
$\Delta _{1}:= \Delta _{0} \cup \mbox{co}_{y_1}F_0  $.
Then take a subfacet $F_1$ in $\Delta _{1}$  and a new vertex $y_2$.
Put
$\Delta _{2}:= \Delta _{1} \cup \mbox{co}_{y_2}F_1  $.
We can continue this process so as to get a sequence of simplicial complexes
$\Delta _{0}, \Delta _{1}, \Delta _{2}, \dots$.
We will refer to any simplicial complex $\Delta _{i}$ ($i\ge 0$) obtained in this way
with the  notation $\Delta + \mbox{($d$-branches)}$ as long as we do not need to mention the sequence 
$F_0,F_1,\dots F_{i-1}$ explicitly.
\newline
 Let $\Delta _1$ and $\Delta_2$ be simplicial complexes on two disjoint vertex sets. 
We define the {\it simplicial join} $\Delta _{1}*\Delta _{2}$ of $\Delta _{1}$
and $\Delta _{2}$ by setting
$$
\Delta _{1}*\Delta _{2}=\{ F \cup G \mid F \in \Delta _{1}, G \in \Delta _{2}\}.
$$
\par\medskip\noindent
\noindent
The following result was proven by the second author of this paper: it follows from \cite{Te}, Theorems 3.2 and 4.2. 
\begin{lemma}\label{terai1}
Let $\Delta$  be a pure $ d$-dimensional strongly connected simplicial
complex.
Then 
$$
{\rm reg } I_{\Delta} \le  \deg K[\Delta]-{\rm  codim }K[\Delta] +1,
$$
and we have
$$
{\rm reg } I_{\Delta} = \deg K[\Delta]-{\rm  codim }K[\Delta] +1 \ge 4
$$
if and only if $\Delta$ can be expressed as
$$ \partial{\Delta}(r)* \Delta(d-r+2)+  
\mbox{\rm ($d$-branches)},
$$
after a suitable change of labeling of the vertices, where $r={\rm reg } I_{\Delta}$.
In particular, in this case  $\Delta$ is Cohen-Macaulay.
\end{lemma}
\par\medskip\noindent
Combining the previous lemma with Theorem \ref{theorem3} we deduce the following result:
\begin{corollary}
Let $\Delta$  be a pure strongly connected 
complex satisfying
$$
{\rm reg } I_{\Delta} = \deg K[\Delta]-{\rm  codim }K[\Delta] +1 \ge 4.
$$
Then  
$I_{\Delta}$ is a set-theoretic complete intersection. 
\end{corollary}
We call a pure $d$-dimensional strongly connected generalized tree a {\it d-tree}. This notion can also be characterized inductively as follows:\newline
(i)  A $d$-simplex is a $d$-tree. \newline
(ii) If $\Delta$ is a $d$-tree, then so is 
$\Delta \cup \mbox{\rm co } _{x_0} F$  
for any  subfacet $F$  of $\Delta$ and any new vertex $x_0$.
\begin{remark}{\rm 
Let $\Delta$  be a pure $d$-dimensional  strongly connected 
complex satisfying
$$
{\rm reg } I_{\Delta} = \deg K[\Delta]-{\rm  codim }K[\Delta] +1\ =2.$$
Then, by Lemma \ref{DiestelFroberg}, $\Delta$  is a generalized tree, hence it is a $d$-tree.   
Thus $I_{\Delta}$ is a set-theoretic complete intersection by Corollary \ref{coro2}.  \newline
On the other hand,
in the case where
${\rm reg } I_{\Delta} = \deg K[\Delta]-{\rm  codim }K[\Delta] +1\ =3$,
we do not know whether $\mbox{ara}\,I_{\Delta} =\mbox{pd}\,K[\Delta]$ still holds. The problem is open 
even in the simple case where $\Delta$ is an $n$-gon: all we know is that the answer is affirmative for $n=3,4$ (where $I_{\Delta}$ is a  complete intersection) and for $n=5,6$ (where $I_{\Delta}$ is  a set-theoretic complete intersection, as was recently proven in \cite{B1}, see also \cite{B0}, Example 2.)}
\end{remark}

\end{document}